\documentclass[12pt]{article}
\usepackage{amsfonts}

\usepackage{mathrsfs}
\usepackage{amscd}
\usepackage{amsmath,amsfonts,amssymb,amscd}
\usepackage{indentfirst,graphics,epsfig,psfrag}
\input{epsf}

\usepackage{amsmath,amssymb,amsthm, amscd, amsfonts, graphicx,ifpdf}
\usepackage{lineno}

\newtheorem{thm}{Theorem}[section]

\newtheorem{lem}{Lemma}[section]
\newtheorem{conj}[thm]{Conjecture}

\def\qed{\nopagebreak\hfill{\rule{4pt}{7pt}}
\medbreak}

\def\pf{\noindent {\it Proof.} }

\title{\bf Hardness results on generalized connectivity\footnote{Supported by NSFC and the ``973" program.}}

\author{
\small Shasha Li, Xueliang Li\\
\small Center for Combinatorics and LPMC-TJKLC\\
\small Nankai University, Tianjin 300071, China.\\
\small  Email: lss@cfc.nankai.edu.cn, lxl@nankai.edu.cn\\
}
\date{}
\begin{document}

\maketitle

\begin{abstract}
Let $G$ be a nontrivial connected graph of order $n$ and let $k$ be
an integer with $2\leq k\leq n$. For a set $S$ of $k$ vertices of
$G$, let $\kappa (S)$ denote the maximum number $\ell$ of
edge-disjoint trees $T_1,T_2,\ldots,T_\ell$ in $G$ such that
$V(T_i)\cap V(T_j)=S$ for every pair $i,j$ of distinct integers with
$1\leq i,j\leq \ell$. A collection $\{T_1,T_2,\ldots,T_\ell\}$ of
trees in $G$ with this property is called an internally disjoint set
of trees connecting $S$. Chartrand et al. generalized the concept of
connectivity as follows: The $k$-$connectivity$, denoted by
$\kappa_k(G)$, of $G$ is defined by
$\kappa_k(G)=$min$\{\kappa(S)\}$, where the minimum is taken over
all $k$-subsets $S$ of $V(G)$. Thus $\kappa_2(G)=\kappa(G)$, where
$\kappa(G)$ is the connectivity of $G$, for which there are
polynomial-time algorithms to solve it.

This paper mainly focus on the complexity of the generalized
connectivity. At first, we obtain that for two fixed positive
integers $k_1$ and $k_2$, given a graph $G$ and a $k_1$-subset $S$
of $V(G)$, the problem of deciding whether $G$ contains $k_2$
internally disjoint trees connecting $S$ can be solved by a
polynomial-time algorithm. Then, we show that when $k_1$ is a fixed
integer of at least $4$, but $k_2$ is not a fixed integer, the
problem turns out to be NP-complete. On the other hand, when $k_2$
is a fixed integer of at least $2$, but $k_1$ is not a fixed
integer, we show that the problem also becomes NP-complete.
Finally we give some open problems.\\[3mm]
\end{abstract}

\section{Introduction}

We follow the terminology and notation of \cite{Bondy} and all
graphs considered here are always simple. The $connectivity$
$\kappa(G)$ of a graph $G$ is defined as the minimum cardinality of
a set $Q$ of vertices of $G$ such that $G-Q$ is disconnected or
trivial. A well-known theorem of Whitney \cite{Whitney} provides an
equivalent definition of connectivity. For each $2$-subset
$S=\{u,v\}$ of vertices of $G$, let $\kappa(S)$ denote the maximum
number of internally disjoint $uv$-paths in $G$. Then
$\kappa(G)=$min$\{\kappa(S)\}$, where the minimum is taken over all
$2$-subsets $S$ of $V(G)$.

In \cite{Chartrand}, the authors generalized the concept of
connectivity. Let $G$ be a nontrivial connected graph of order $n$
and let $k$ be an integer with $2\leq k\leq n$. For a set $S$ of $k$
vertices of $G$, let $\kappa (S)$ denote the maximum number $\ell$
of edge-disjoint trees $T_1,T_2,\ldots,T_\ell$ in $G$ such that
$V(T_i)\cap V(T_j)=S$ for every pair $i,j$ of distinct integers with
$1\leq i,j\leq \ell$ (Note that the trees are vertex-disjoint in
$G\backslash S$). A collection $\{T_1,T_2,\ldots,T_\ell \}$ of trees
in $G$ with this property is called an {\it internally disjoint set
of trees connecting $S$}. The $k$-$connectivity$, denoted by
$\kappa_k(G)$, of $G$ is then defined by
$\kappa_k(G)=$min$\{\kappa(S)\}$, where the minimum is taken over
all $k$-subsets $S$ of $V(G)$. Thus, $\kappa_2(G)=\kappa(G)$.

In \cite{LLZ}, we focused on the investigation of $\kappa_3(G)$ and
mainly studied the relationship between the $2$-connectivity and the
$3$-connectivity of a graph. We gave sharp upper and lower bounds of
$\kappa_3(G)$ for general graphs $G$, and constructed two kinds of
graphs which attain the upper and lower bound, respectively. We also
showed that if $G$ is a connected planar graph, then $\kappa(G)-1
\leq \kappa_3(G)\leq \kappa(G)$, and gave some classes of graphs
which attain the bounds. Moreover, we studied algorithmic aspects
for $\kappa_3(G)$ and gave an algorithm to determine $\kappa_3(G)$
for general graph $G$. This algorithm runs in a polynomial time for
graphs with a fixed value of connectivity, which implies that the
problem of determining $\kappa_3(G)$ for graphs with a small minimum
degree or connectivity can be solved in polynomial time, in
particular, the problem whether $\kappa(G)=\kappa_3(G)$ for a planar
graph $G$ can be solved in polynomial time.

In this paper, we will turn to the complexity of the generalized
connectivity. At first, by generalizing the algorithm of \cite{LLZ},
we obtain that for two fixed positive integers $k_1$ and $k_2$,
given a graph $G$ and a $k_1$-subset $S$ of $V(G)$, the problem of
deciding whether $G$ contains $k_2$ internally disjoint trees
connecting $S$ can be solved by a polynomial-time algorithm. Then,
we show that when $k_1$ is a fixed integer of at least $4$, but
$k_2$ is not a fixed integer, the problem turns out to be
NP-complete.

\begin{thm}\label{thm1}
For any fixed integer $k_1\geq 4$, given a graph $G$, a $k_1$-subset
$S$ of $V(G)$ and an integer $2\leq k_2\leq n-1$, deciding whether
there are $k_2$ internally disjoint trees connecting $S$, namely
deciding whether $\kappa (S)\geq k_2$, is NP-complete.
\end{thm}

On the other hand, when $k_2$ is a fixed integer of at least $2$,
but $k_1$ is not a fixed integer, we show that the problem also
becomes NP-complete.

\begin{thm}\label{thm2}
For any fixed integer $k\geq 2$, given a graph $G$ and a subset $S$
of $V(G)$, deciding whether there are $k$ internally disjoint trees
connecting $S$, namely deciding whether $\kappa (S)\geq k$, is
NP-complete.
\end{thm}

The rest of this paper is organized as follows. The next section
simply generalizes the algorithm of \cite{LLZ} and makes some
preparations. Then Sections $3$ and $4$ prove Theorem \ref{thm1} and
Theorem \ref{thm2}, respectively. The final section, Section $5$,
contains some open problems.

\section{Preliminaries}

At first, we introduce the following result of \cite{LLZ}.

\begin{lem}\label{lem1}
Given a fixed positive integer $k$, for any graph $G$ the problem of
deciding whether $G$ contains $k$ internally disjoint trees
connecting $\{v_1,v_2,v_3\}$ can be solved by a polynomial-time
algorithm, where $v_1,v_2,v_3$ are any three vertices of $V(G)$.
\end{lem}

We first show that the trees we really want has only two types. Then
we prove that if there are $k$ internally disjoint trees connecting
$\{v_1,v_2,v_3\}$, then the union of the $k$ trees has at most
$f(k)n^k$ types, where $f(k)$ is a function on $k$. For every $i\in
[f(k)n^k]$, we can convert into a $k'$-linkage problem the problem
of deciding whether $G$ contains a union of $k$ trees having type
$i$. Since the $k'$-linkage problem has a polynomial-time algorithm
to solve it, which has a running time $O(n^3)$, see
\cite{Robertson}, and $k$ is a fixed integer, we finally obtain that
the problem of deciding whether $\kappa{\{v_1,v_2,v_3\}}\geq k$ can
be solved by a polynomial-time algorithm. We refer the readers to
\cite{LLZ} for details.

By the similar method, we can also show that
given a fixed positive integer $k$,
for any graph $G$ the problem of deciding whether $G$ contains $k$
internally disjoint trees connecting $\{v_1,v_2,v_3,v_4\}$
can be solved by a polynomial-time algorithm, where $v_1,v_2,v_3,v_4$
are any four vertices of $V(G)$.

\begin{figure}[h,t]
\begin{center}
\input{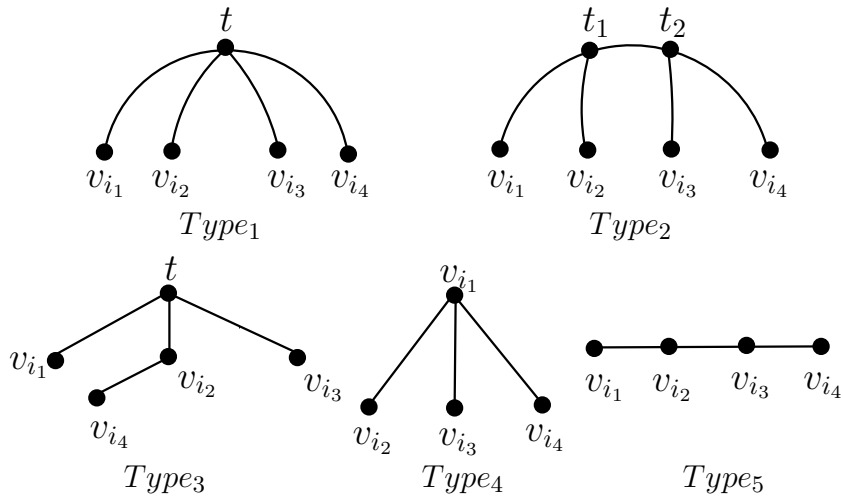}
\caption{Five types of trees we really want, where
$\{v_{i_1},v_{i_2},v_{i_3},v_{i_4}\}=\{v_1,v_2,v_3,v_4\}.$}
\label{fig1}
\end{center}
\end{figure}

Since for the trees $T$ connecting $\{v_1,v_2,v_3,v_4\}$, we only
need $T$ belonging to one of the five types in Figure \ref{fig1},
then if there are $k$ internally disjoint trees connecting
$\{v_1,v_2,v_3,v_4\}$, consider the union of the $k$ trees and it is
not hard to obtain that the number of types is at most $f(k)n^{2k}$,
where $f(k)$ is a function on $k$ and $f(k)n^{2k}$ is only a rough
upper bound. Then for every $i\in [f(k)n^{2k}]$, we can convert into
a $k'$-linkage problem the problem of deciding whether $G$ contains
a union of $k$ trees having type $i$. Since the $k'$-linkage problem
has a polynomial-time algorithm and $k$ is a fixed integer, we
obtain that the problem of deciding whether
$\kappa{\{v_1,v_2,v_3,v_4\}}\geq k$ can be solved by a
polynomial-time algorithm.

Now, for two fixed positive integers $k_1$ and $k_2$, if we replace
the set $\{v_1,v_2,v_3,v_4\}$ with a $k_1$-subset $S$ of $V(G)$ and
replace $k$ with $k_2$, the problem can still be solved by a
polynomial-time algorithm. The method is similar.

Since for the trees $T$ connecting the $k_1$-subset $S$ of $V(G)$,
the number of types of $T$ we really want is at most $f_1(k_1)$,
where $f_1(k_1)$ is a function on $k_1$, then if there are $k_2$
internally disjoint trees connecting $S$, consider the union of the
$k_2$ trees and it is not hard to obtain that the number of types is
at most $f_2(k_1,k_2)n^{(k_1-2)k_2}$, where $f_2(k_1,k_2)$ is a
function on $k_1$ and $k_2$ and $f_2(k_1,k_2)n^{(k_1-2)k_2}$ is only
a rough upper bound. Next, by the same way, for every $i\in
[f_2(k_1,k_2)n^{(k_1-2)k_2}]$, convert into a $k'$-linkage problem
the problem of deciding whether $G$ contains a union of $k_2$ trees
having type $i$ and a polynomial-time algorithm is then obtained.

\begin{lem}\label{lem2}
For two fixed positive integer $k_1$ and $k_2$,
given a graph $G$ and a $k_1$-subset $S$ of
$V(G)$, the problem of deciding whether $G$ contains $k_2$
internally disjoint trees connecting $S$
can be solved by a polynomial-time algorithm.
\end{lem}

Note that Lemma \ref{lem2} is a generalization of Lemma \ref{lem1}.
When $k_1=3$ and $k_2=k$, Lemma \ref{lem2}
is exactly Lemma \ref{lem1}.

Before proceeding, we recall the following two basic NP-complete
problems.

\noindent {\bf $3$-DIMENSIONAL MATCHING ($3$-DM)}

\noindent Given three sets $U$, $V$, and $W$ of equal cardinality,
and a subset $T$ of $U\times V\times W$, decide whether there is a
subset $M$ of $T$ with $|M|=|U|$ such that whenever $(u,v,w)$ and
$(u',v',w')$ are distinct triples in $M$, $u\neq u'$, $v\neq v'$,
and $w\neq w'$ ?

\noindent {\bf BOOLEAN $3$-SATISFIABILITY ($3$-SAT)}

\noindent Given a boolean formula $\phi$ in conjunctive normal form
with three literals per clause, decide whether $\phi$ is satisfiable
?

\section{Proof of Theorem \ref{thm1}}

For the problem in Lemma \ref{lem2}, when $k_1=4$ and $k_2$ is not a
fixed integer, we denote this case by Problem $1$.

\noindent {\bf Problem $1$.} Given a graph $G$, a $4$-subset $S$ of
$V(G)$ and an integer $2\leq k\leq n-1$, decide whether there are
$k$ internally disjoint trees connecting $S$, namely decide whether
$\kappa (S)\geq k$ ?

At first, we will show that Problem $1$ is NP-complete by reducing
$3$-DM to it, as follows.

\begin{lem}\label{lem3}
Given a graph $G$, a $4$-subset $S$ of $V(G)$ and an integer $2\leq
k\leq n-1$, deciding whether there are $k$ internally disjoint trees
connecting $S$, namely deciding whether $\kappa (S)\geq k$, is
NP-complete.
\end{lem}
\pf It is clear that Problem $1$ is in NP. So it will suffice to show that
$3$-DM is polynomially reducible to this problem.

Given three sets of equal cardinality, denoted by $U=\{u_1,u_2,\ldots,u_n\}$,
$V=\{v_1,v_2,\ldots,v_n\}$ and $W=\{w_1,w_2,\ldots,w_n\}$, and a subset
$T=\{T_1,T_2,\ldots,T_m\}$ of $U\times V\times W$, we will construct a graph
$G'$, a $4$-subset $S$ of $V(G')$ and an integer $k\leq |V(G')|-1$
such that there are $k$ internally disjoint trees connecting $S$
in $G'$ if and only if there is a subset $M$ of $T$ with $|M|=|U|=n$
such that whenever $(u_i,v_j,w_k)$ and $(u_{i'},v_{j'},w_{k'})$
are distinct triples in $M$, $u_i\neq u_{i'}$, $v_j\neq v_{j'}$ and
$w_k\neq w_{k'}$.

We define $G'$ as follows:

\begin{align*}
V(G')&=\{\hat{u},\hat{v},\hat{w},\hat{t}\}\cup
\{u_i:1\leq i\leq n\}\cup \{v_i:1\leq i\leq n\}\\
&\cup \{w_i:1\leq i\leq n\}\cup \{t_i:1\leq i\leq m\}
\cup \{a_i:1\leq i\leq m-n\};\\
E(G')&=\{\hat{u}u_i:1\leq i\leq n\}\cup \{\hat{v}v_i:1\leq i\leq n\}
\cup \{\hat{w}w_i:1\leq i\leq n\}\\
&\cup\{\hat{t}t_i:1\leq i\leq m\}\cup \{\hat{u}a_i:1\leq i\leq m-n\}
\cup \{\hat{v}a_i:1\leq i\leq m-n\}\\
&\cup \{\hat{w}a_i:1\leq i\leq m-n\}\cup \{t_ia_j: 1\leq i\leq m, 1\leq j\leq m-n\}\\
&\cup \{t_iu_j:u_j\in T_i\} \cup \{t_iv_j:v_j\in T_i\}
\cup \{t_iw_j:w_j\in T_i\}.
\end{align*}

Then let $S=\{\hat{u},\hat{v},\hat{w},\hat{t}\}$ and $k=m$.

Suppose that there is a subset $M$ of $T$ with $|M|=|U|=n$
such that whenever $(u_i,v_j,w_k)$ and $(u_{i'},v_{j'},w_{k'})$
are distinct triples in $M$, $u_i\neq u_{i'}$, $v_j\neq v_{j'}$ and
$w_k\neq w_{k'}$. Then for every $T_i\in M$, we can construct a tree
whose vertex set consists of $S$, $t_i$ and three vertices corresponding
to three elements in $T_i$. For each $T-i\notin M$,
$G[t_i,a_j,\hat{u},\hat{v},\hat{w},\hat{t}]$ is a tree connecting
$S$, for some $1\leq j\leq m-n$. So we can easily find out
$k$ internally disjoint trees connecting $S$ in $G'$.

Now suppose that there are $k=m$ internally disjoint trees
connecting $S$ in $G'$. Since $\hat{u},\hat{v},\hat{w}$ and
$\hat{t}$ all have degree $m$, then among the $m$ trees, there are
$n$ trees, each of which contains the vertices in $S$, a vertex from
$\{t_i:1\leq i\leq m\}$, a vertex from $\{u_i:1\leq i\leq n\}$, a
vertex from $\{v_i:1\leq i\leq n\}$ and a vertex from $\{w_i:1\leq
i\leq n\}$ and can not contain any other vertex. Since the $n$ trees
are internally disjoint, it can be easily checked that $n$ $3$-sets
$T_i\in U\times V\times W$ corresponding to $n$ vertices $t_i$ in
the $n$ trees form a subset $M$ of $T$ with $|M|=|U|=n$ such that
whenever $(u_i,v_j,w_k)$ and $(u_{i'},v_{j'},w_{k'})$ are distinct
triples in $M$, $u_i\neq u_{i'}$, $v_j\neq v_{j'}$ and $w_k\neq
w_{k'}$. The proof is complete.\qed

Now we show that for a fixed integer $k_1\geq 5$, in Problem $1$
replacing the $4$-subset of $V(G)$ with a $k_1$-subset of $V(G)$,
the problem is still NP-complete, which can easily be proved by
reducing Problem $1$ to it.

\begin{lem}\label{lem4}
For any fixed integer $k_1\geq 5$, given a graph $G$, a $k_1$-subset
$S$ of $V(G)$ and an integer $2\leq k_2\leq n-1$, deciding whether
there are $k_2$ internally disjoint trees connecting $S$, namely
deciding whether $\kappa (S)\geq k_2$, is NP-complete.
\end{lem}
\pf Clearly, the problem is in NP. We will prove that Problem $1$ is
polynomially reducible to it.

For any given graph $G$, a $4$-subset $S=\{v_1,v_2,v_3,v_4\}$ of
$V(G)$ and an integer $2\leq k\leq n-1$, we construct a new graph
$G'=(V',E')$ and a $k_1$-subset $S'$ of $V(G')$ and let $k_2=k$ be
such that there are $k_2=k$ internally disjoint trees connecting
$S'$ in $G'$ if and only if there are $k$ internally disjoint trees
connecting $S$ in $G$.

We construct $G'=(V',E')$ by adding $k_1-4$ new vertices
$\{\hat{a}^1,\hat{a}^2,\ldots,\hat{a}^{k_1-4}\}$ to $G$ and for
every $i\leq k_1-4$, adding $k_2$ internally disjoint
$\hat{a}^iv_1$-paths $\{\hat{a}^ia^i_jv_1:1\leq j\leq k_2\}$ of
length two, where $a^i_j$ is also a new vertex and if $i_1\neq i_2$,
$a^{i_1}_{j_1}\neq a^{i_2}_{j_2} $. Then let
$S'=\{v_1,v_2,v_3,v_4,\hat{a}^1,\hat{a}^2,\ldots,\hat{a}^{k_1-4}\}$.
It is not hard to check that $\kappa_{G'} (S')\geq k_2=k$ if and
only if $\kappa_G (S)\geq k$. The proof is complete.\qed

Combining Lemma \ref{lem3} with Lemma \ref{lem4},
we obtain Theorem \ref{thm1}, namely, we complete the proof
of Theorem \ref{thm1}.

\section{Proof of Theorem \ref{thm2}}

For the problem in Lemma \ref{lem2}, when $k_2=2$ and $k_1$ is not a
fixed integer, we denote this case by Problem $2$.

\noindent {\bf Problem $2$.} Given a graph $G$ and a subset $S$ of
$V(G)$, decide whether there are two internally disjoint trees
connecting $S$, namely decide whether $\kappa (S)\geq 2$ ?

Firstly, the following lemma will prove that Problem $2$ is
NP-complete by reducing $3$-SAT to it.

\begin{lem}\label{lem5}
Given a graph $G$ and a subset $S$ of $V(G)$, deciding whether there
are two internally disjoint trees connecting $S$, namely deciding
whether $\kappa (S)\geq 2$, is NP-complete.
\end{lem}
\pf Clearly, Problem $2$ is in NP. So it will suffice to show that
$3$-SAT is polynomially reducible to this problem.

Given a $3$-CNF formula $\phi=\bigwedge_{i=1}^m c_i$ over variables
$x_1,x_2,\ldots,x_n$, we construct a graph $G_{\phi}$ and a subset
$S$ of $V(G_{\phi})$ such that there are two internally disjoint
trees connecting $S$ if and only if $\phi$ is satisfiable.

We define $G_{\phi}$ as follows:
\begin{align*}
V(G_{\phi})&=\{\hat{x_i}: 1\leq i\leq n\}\cup \{x_i: 1\leq i\leq
n\}\cup \{\ \bar{x_i}: 1\leq i\leq n\} \\
&\cup \{c_i: 1\leq i\leq m\}\cup \{a\};\\
E(G_{\phi})&=\{\hat{x_i}x_i: 1\leq i\leq n\}\cup
\{\hat{x_i}\bar{x_i}: 1\leq i\leq n\}\\
&\cup \{x_ic_j: x_i\in c_j\}\cup \{\bar{x_i}c_j: \bar{x_i}\in
c_j\}\\
&\cup \{x_1x_i: 2\leq i\leq n\}\cup \{x_1\bar{x_i}: 2\leq i\leq
n\}\cup \{\bar{x_1}x_i: 2\leq i\leq n\}\cup \{\bar{x_1}\bar{x_i}:
2\leq
i\leq n\}\\
&\cup \{ax_i: 1\leq i\leq n\}\cup \{a\bar{x_i}: 1\leq i\leq n\}\cup
\{ac_i: 1\leq i\leq m\},
\end{align*}
where the notation $x_i\in c_j$($\bar{x_i}\in c_j$) signifies that
$x_i$($\bar{x_i}$) is a literal of the clause $c_j$. Then let
$S=\{\hat{x_i}: 1\leq i\leq n\}\cup \{c_i: 1\leq i \leq m\}$.

Suppose that there is a true assignment $t$ satisfying $\phi$. Then
for every clause $c_i$($1\leq i\leq m$), there must exist a literal
$x_j\in c_i$ such that $t(x_j)=1$ or $\bar{x_j}\in c_i$ such that
$t(x_j)=0$, for some $1\leq j\leq m$. For such literals $x_j$ or
$\bar{x_j}$, let $T_1$ be a graph such that $E(T_1)=\{c_ix_j$ (or
$c_i\bar{x_j}): 1\leq i\leq m\}$. Obviously, at most one of the two
vertices $x_j$ and $\bar{x_j}$ exists in $V(T_1)$. If neither $x_j$
nor $\bar{x_j}$ is in $V(T_1)$, we can add any one of them to
$V(T_1)$. Now, if $x_1\in V(T_1)$, add $x_1x_i$(if $x_i\in V(T_1)$)
or $x_1\bar{x_i}$(if $\bar{x_i}\in V(T_1)$) to $E(T_1)$, for $2\leq
i\leq n$. Otherwise, add $\bar{x_1}x_i$(if $x_i\in V(T_1)$) or
$\bar{x_1}\bar{x_i}$(if $\bar{x_i}\in V(T_1)$) to $E(T_1)$. Finally,
add edges $x_i\hat{x_i}$(if $x_i\in V(T_1)$) or
$\bar{x_i}\hat{x_i}$(if $\bar{x_i}\in V(T_1)$) to $E(T_1)$, for
$1\leq i\leq n$. Now it is easy to check that $T_1$ is a tree
connecting $S$. Then let $T_2$ be a tree containing $ac_i$ for
$1\leq i \leq m$, $ax_j$ and $x_j\hat{x_j}$(if $\bar{x_j}\in
V(T_1)$) or $a\bar{x_j}$ and $\bar{x_j}\hat{x_j}$(if $x_j\in
V(T_1)$) for $1\leq j \leq n$. $T_1$ and $T_2$ are two internally
disjoint trees connecting $S$.

Now suppose that there are two internally disjoint trees $T_1,T_2$
connecting $S$. Since $a\notin S$, only one tree can contain the
vertex $a$. Without loss of generality, assume that $a\notin
V(T_1)$. Since for every $1\leq i\leq n$, $\hat{x_i}\in S$ has
degree two, $V(T_1)$ must contain one and only one of its two
neighbors $x_i$ and $\bar{x_i}$. Then let the value of a variable
$x_i$ be $1$ if its corresponding vertex $x_i$ is contained in
$V(T_1)$. Otherwise let the value be $0$. Moreover, because $a\notin
V(T_1)$, for every $c_i$($1\leq i\leq m$), there must exist some
vertex $x_j\in V(T_1)$ such that $c_ix_j\in E(T_1)$ or $\bar{x_j}\in
V(T_1)$ such that $c_i\bar{x_j}\in E(T_1)$. So, $\phi$ is obviously
satisfiable by the above true assignment. The proof is complete.\qed

Now we show that for a fixed integer $k\geq 3$, in Problem $2$ if we
want to decide whether there are $k$ internally disjoint trees
connecting $S$ rather than two, the problem is still NP-complete,
which can easily be proved by reducing Problem $2$ to it.

\begin{lem}\label{lem6}
For any fixed integer $k\geq 3$, given a graph $G$ and a subset $S$
of $V(G)$, deciding whether there are $k$ internally disjoint trees
connecting $S$, namely deciding whether $\kappa (S)\geq k$, is
NP-complete.
\end{lem}
\pf Clearly, the problem is in NP. We will show that Problem $2$ is
polynomially reducible to this problem.

Note that $k$ is an fixed integer of at least $3$. For any given
graph $G$ and a subset $S$ of $V(G)$, we construct a graph
$G'=(V',E')$ by adding $k-2$ new vertices to $G$ and joining every
new vertex to all vertices in $S$. Then let $S'$ be a subset of
$V(G')$ such that $S'=S$.

If $\kappa_G(S)\geq 2$, it is clear that $\kappa_{G'}(S')\geq k$.

Suppose that there are $k$ internally disjoint trees connecting $S'$
in $G'$, namely $\kappa_{G'}(S')\geq k$. Since there are only $k-2$
new vertices, at least two trees can not contain any new vertex,
which means the two trees are actually two internally disjoint trees
connecting $S'=S$ in $G$. The proof is complete.\qed

Combining Lemma \ref{lem5} with Lemma \ref{lem6},
we obtain Theorem \ref{thm2}, namely, we complete the proof
of Theorem \ref{thm2}.

\section{Open problems}

As Theorem \ref{thm1}, we only show that for any fixed integer
$k_1\geq 4$, given a graph $G$, a $k_1$-subset $S$ of $V(G)$ and an
integer $2\leq k_2\leq n-1$, deciding whether $\kappa (S)\geq k_2$
is NP-complete, while for $k_1=3$, the complexity is not known.
However, we tend to believe that it is NP-complete.

\begin{conj}\label{conj1}
Given a graph $G$, a $3$-subset $S$ of $V(G)$ and an integer $2\leq
k\leq n-1$, deciding whether there are $k$ internally disjoint trees
connecting $S$, namely deciding whether $\kappa (S)\geq k$ is
NP-complete.
\end{conj}

By Lemma \ref{lem1}, we know that given a fixed positive integer
$k$, for any graph $G$ and a $3$-subset $S$ of $V(G)$ the problem of
deciding whether $\kappa(S)\geq k$ can be solved by a
polynomial-time algorithm. Moreover, by the definition
$\kappa_3(G)=min\{\kappa(S)\}$, where the minimum is taken over all
$3$-subsets $S$ of $V(G)$, we therefore obtain that the problem of
deciding whether $\kappa_3(G)\geq k$ can also be solved by a
polynomial-time algorithm \cite{LLZ}.

Similarly, since we know that given two fixed integers $k_1\geq 4$
and $k_2$, for any graph $G$ and a $k_1$-subset $S$ of $V(G)$ the
problem of deciding whether $\kappa(S)\geq k_2$ can be solved by a
polynomial-time algorithm and $\kappa_{k_1}(G)=min\{\kappa(S)\}$,
where the minimum is taken over all $k_1$-subsets $S$ of $V(G)$, we
can also obtain that the problem of deciding whether
$\kappa_{k_1}(G)\geq k_2$ can be solved by a polynomial-time
algorithm.

However, if $k_2$ is not a fixed positive integer, the complexity of
the problem is still not known, including the case of $k_1=3$. We
conjecture that it could be NP-complete, as follows.

\begin{conj}\label{conj2}
For a fixed integer $k_1\geq 3$, given a graph $G$ and an integer
$2\leq k\leq n-1$, the problem of deciding whether
$\kappa_{k_1}(G)\geq k$ is NP-complete.
\end{conj}

\end{document}